\documentclass{amsart}
\usepackage{latexsym}
\usepackage{amssymb}
\usepackage{amsmath}
\begin{document}

%%%%%%%%%%%%%%%%%%%%%%Definitions%%%%%%%%%%%%%%%%%%%%%%%%%%%%%%%%%%%%%%%%%%%

\newtheorem{theorem}{Theorem}
\newtheorem{problem}{Problem}
\newtheorem{definition}{Definition}
\newtheorem{lemma}{Lemma}
\newtheorem{proposition}{Proposition}
\newtheorem{corollary}{Corollary}
\newtheorem{example}{Example}
\newtheorem{conjecture}{Conjecture}
\newtheorem{algorithm}{Algorithm}
\newtheorem{exercise}{Exercise}
\newtheorem{xample}{Example}
\newtheorem{remarkk}{Remark}

\newcommand{\be}{\begin{equation}}
\newcommand{\ee}{\end{equation}}
\newcommand{\bea}{\begin{eqnarray}}
\newcommand{\eea}{\end{eqnarray}}
\newcommand{\beq}[1]{\begin{equation}\label{#1}}
\newcommand{\eeq}{\end{equation}}
\newcommand{\beqn}[1]{\begin{eqnarray}\label{#1}}
\newcommand{\eeqn}{\end{eqnarray}}
\newcommand{\beaa}{\begin{eqnarray*}}
\newcommand{\eeaa}{\end{eqnarray*}}
\newcommand{\req}[1]{(\ref{#1})}

\newcommand{\lip}{\langle}
\newcommand{\rip}{\rangle}
\newcommand{\uu}{\underline}
\newcommand{\oo}{\overline}
\newcommand{\La}{\Lambda}
\newcommand{\la}{\lambda}
\newcommand{\eps}{\varepsilon}
\newcommand{\om}{\omega}
\newcommand{\Om}{\Omega}
\newcommand{\ga}{\gamma}
\newcommand{\ka}{\kappa}
\newcommand{\rrr}{{\Bigr)}}
\newcommand{\qqq}{{\Bigl\|}}

\newcommand{\dint}{\displaystyle\int}
\newcommand{\dsum}{\displaystyle\sum}
\newcommand{\dfr}{\displaystyle\frac}
\newcommand{\bige}{\mbox{\Large\it e}}
\newcommand{\integers}{{\Bbb Z}}
\newcommand{\rationals}{{\Bbb Q}}
\newcommand{\reals}{{\rm I\!R}}
\newcommand{\realsd}{\reals^d}
\newcommand{\realsn}{\reals^n}
\newcommand{\NN}{{\rm I\!N}}
\newcommand{\DD}{{\rm I\!D}}
\newcommand{\degree}{{\scriptscriptstyle \circ }}
\newcommand{\dfn}{\stackrel{\triangle}{=}}
\def\complex{\mathop{\raise .45ex\hbox{${\bf\scriptstyle{|}}$}
     \kern -0.40em {\rm \textstyle{C}}}\nolimits}
\def\hilbert{\mathop{\raise .21ex\hbox{$\bigcirc$}}\kern -1.005em {\rm\textstyle{H}}} %Hilbert space
\newcommand{\RAISE}{{\:\raisebox{.6ex}{$\scriptstyle{>}$}\raisebox{-.3ex}
           {$\scriptstyle{\!\!\!\!\!<}\:$}}} % >< one above each other

\newcommand{\hh}{{\:\raisebox{1.8ex}{$\scriptstyle{\degree}$}\raisebox{.0ex}
           {$\textstyle{\!\!\!\! H}$}}}

\newcommand{\OO}{\won}
\newcommand{\calA}{{\mathcal A}}
\newcommand{\calB}{{\mathcal B}}
\newcommand{\calC}{{\cal C}}
\newcommand{\calD}{{\mathcal D}}
\newcommand{\calE}{{\cal E}}
\newcommand{\calF}{{\mathcal F}}
\newcommand{\calG}{{\cal G}}
\newcommand{\calH}{{\cal H}}
\newcommand{\calK}{{\cal K}}
\newcommand{\calL}{{\mathcal L}}
\newcommand{\calM}{{\mathcal M}}
\newcommand{\calO}{{\cal O}}
\newcommand{\calP}{{\cal P}}
\newcommand{\calS}{{\mathcal S}}
\newcommand{\calT}{{\mathcal T}}
\newcommand{\calU}{{\mathcal U}}
\newcommand{\calX}{{\cal X}}
\newcommand{\calZ}{{\mathcal Z}}
\newcommand{\calXX}{{\cal X\mbox{\raisebox{.3ex}{$\!\!\!\!\!-$}}}}
\newcommand{\calXXX}{{\cal X\!\!\!\!\!-}}
\newcommand{\gi}{{\raisebox{.0ex}{$\scriptscriptstyle{\cal X}$}
\raisebox{.1ex} {$\scriptstyle{\!\!\!\!-}\:$}}}
\newcommand{\intsim}{\int_0^1\!\!\!\!\!\!\!\!\!\sim}
\newcommand{\intsimt}{\int_0^t\!\!\!\!\!\!\!\!\!\sim}
\newcommand{\pp}{{\partial}}
\newcommand{\al}{{\alpha}}
\newcommand{\sB}{{\cal B}}
\newcommand{\sL}{{\cal L}}
\newcommand{\sF}{{\cal F}}
\newcommand{\sE}{{\cal E}}
\newcommand{\sX}{{\cal X}}
\newcommand{\R}{{\rm I\!R}}
\renewcommand{\L}{{\rm I\!L}}
\newcommand{\vp}{\varphi}
\newcommand{\N}{{\rm I\!N}}
\def\ooo{\lip}
\def\ccc{\rip}
\newcommand{\ot}{\hat\otimes}
\newcommand{\rP}{{\Bbb P}}
\newcommand{\bfcdot}{{\mbox{\boldmath$\cdot$}}}

\renewcommand{\varrho}{{\ell}}
\newcommand{\dett}{{\textstyle{\det_2}}}
\newcommand{\sign}{{\mbox{\rm sign}}}
\newcommand{\TE}{{\rm TE}}
\newcommand{\TA}{{\rm TA}}
\newcommand{\E}{{\rm E\,}}
\newcommand{\won}{{\mbox{\bf 1}}}
\newcommand{\Lebn}{{\rm Leb}_n}
\newcommand{\Prob}{{\rm Prob\,}}
\newcommand{\sinc}{{\rm sinc\,}}
\newcommand{\ctg}{{\rm ctg\,}}
\newcommand{\loc}{{\rm loc}}
\newcommand{\trace}{{\,\,\rm trace\,\,}}
\newcommand{\Dom}{{\rm Dom}}
\newcommand{\ifff}{\mbox{\ if and only if\ }}
\newcommand{\nproof}{\noindent {\bf Proof:\ }}
\newcommand{\remark}{\noindent {\bf Remark:\ }}
\newcommand{\remarks}{\noindent {\bf Remarks:\ }}
\newcommand{\note}{\noindent {\bf Note:\ }}

\newcommand{\boldx}{{\bf x}}
\newcommand{\boldX}{{\bf X}}
\newcommand{\boldy}{{\bf y}}
\newcommand{\boldR}{{\bf R}}
\newcommand{\uux}{\uu{x}}
\newcommand{\uuY}{\uu{Y}}

\newcommand{\limn}{\lim_{n \rightarrow \infty}}
\newcommand{\limN}{\lim_{N \rightarrow \infty}}
\newcommand{\limr}{\lim_{r \rightarrow \infty}}
\newcommand{\limd}{\lim_{\delta \rightarrow \infty}}
\newcommand{\limM}{\lim_{M \rightarrow \infty}}
\newcommand{\limsupn}{\limsup_{n \rightarrow \infty}}

\newcommand{\ra}{ \rightarrow }

\newcommand{\ARROW}[1]
  {\begin{array}[t]{c}  \longrightarrow \\[-0.2cm] \textstyle{#1} \end{array} }

\newcommand{\AR}
 {\begin{array}[t]{c}
  \longrightarrow \\[-0.3cm]
  \scriptstyle {n\rightarrow \infty}
  \end{array}}

\newcommand{\pile}[2]
  {\left( \begin{array}{c}  {#1}\\[-0.2cm] {#2} \end{array} \right) }

\newcommand{\floor}[1]{\left\lfloor #1 \right\rfloor}

%for doing boldface subscripts etc., e.g. $G_{\mmbox{\boldx}}$
\newcommand{\mmbox}[1]{\mbox{\scriptsize{#1}}}

%fraction with round brackets
\newcommand{\ffrac}[2]
  {\left( \frac{#1}{#2} \right)}

\newcommand{\one}{\frac{1}{n}\:}
\newcommand{\half}{\frac{1}{2}\:}

\def\le{\leq}
\def\ge{\geq}
\def\lt{<}
\def\gt{>}

%qed
\def\squarebox#1{\hbox to #1{\hfill\vbox to #1{\vfill}}}
\newcommand{\nqed}{\hspace*{\fill}
          \vbox{\hrule\hbox{\vrule\squarebox{.667em}\vrule}\hrule}\bigskip}
%%%%%%%%%%%%%%%%%%%%%%%%%%%%%%%%%%%%%%%%%%%%%%%%%%%%%%%%%%%%%%%%%%%%%%%%%%%%

\title{Probabilistic Solution of the American Options}

\author{AL\.{I}  S\"uleyman  \"Ust\"unel}
%\date{ }
\maketitle
\begin{abstract}
The existence and uniqueness of probabilistic solutions of variational
inequalities for the general American options are proved under the
hypothesis of hypoellipticity of the infinitesimal generator of the
underlying diffusion process which represents the  risky assets of the
stock market with which the option is created. The main tool is an extension of the It\^o formula which
is valid for the tempered distributions on $\R^d$ and for  nondegenerate
It\^o processes in the sense of the Malliavin calculus.
\end{abstract}

%\tableofcontents
\section{Introduction}
The difficulty to justify the validity of the probabilistic solutions
of the American options is well-known. This is in fact due to the lack
of regularity of the classical solutions of the variational
inequalities (cf.\cite{B-L}) which are satisfied by the value function which
characterizes the Snell envelope (cf. \cite{Lam} for a recent survey
about this subject). In particular the value function is not twice
differentiable hence  the It\^o formula is not applicable to apply the
usual probabilistic techniques. In the case of Black and Scholes model, there are
some results   using extensions of the It\^o formula for the Brownian
motion, which, however, are of limited utility for more general cases.

In this note we give hopefully more general results in the
sense that the option is constructed by the assets which obey to a
general, finite dimensional stochastic differential equation with
deterministic coefficients, i.e.,  a diffusion process. The basic hypothesis used is the
nondegeneracy of  this diffusion in the sense of the Malliavin
calculus  (cf. \cite{Mal}): recall that an $\R^d$-valued  random variable $F=(F_1,\ldots,F_d)$, 
defined on a Wiener space is called nondegenerate
(cf.\cite{Mal,ASU,ASU-1})  if it is infinitely
Sobolev differentiable with respect to the Wiener measure and if the
determinant of the inverse of the matrix $((\nabla F_i,\nabla
F_j)_H\,:i,\,j\leq d)$, where $\nabla$ denotes the Sobolev derivative
on the Wiener space, is in all the $L^p$-spaces w.r. to the Wiener
measure. In this case, the mapping $f\to f\circ F$, defined from  the smooth functions
on $\R^d$ to the space of smooth functions on the Wiener space extends
continuously to a linear  mapping, denoted as $T\to
T(F),\,T\in\calS'(\R^d)$,  from the tempered distributions
$\calS'(\R^d)$ to the space of Meyer distributions on the Wiener space
(cf.\cite{Mal,ASU,ASU-1}). Similarly, if $F$ is replaced with  an It\^o
process whose components satisfy similar regularity  properties, we obtain an
It\^o formula for $T(F_t)-T(F_s)$, $0<s\leq t$,  where the stochastic integral should
be treated as a distribution-valued Gaussian divergence and the
absolutely continuous term is a Bochner integral concentrated in some
negatively indexed Sobolev space. Moreover, if this
latter term  is a  positive distribution, then the resulting integral is a Radon measure on
the Wiener space due to a well-known result about the positive Meyer
distributions on the Wiener space (cf.\cite{ASU-0,ASU-00,ASU,ASU-1}). 

Having summarized the  technical tools  that we use, let us
explain now the main results of the paper: for the uniqueness result  we treat two different
situations; namely the first one where the coefficients are time
dependent and  the variational inequality is
interpreted as an evolutionary variational inequality in
$\calS'(\R^d)$. The second one concerns the case where the
coefficients are time-independent and we interpret it as an inequality
in the space $\calD'(0,T)\otimes \calS'(\R^d)$ with a boundary
condition,  which is of course
more general than the first one. In both cases the operators are
supposed only to be  hypoelliptic;  a hypothesis which is far more general
than the ellipticity hypothesis used in \cite{B-L}. The homogeneity in
time permits us more generality since, in this case the time-component
regularization by the 
mollifiers  of the  solution candidates preserve their property of
being negative distributions, hence measures. The existence is studied
in the last section using the similar techniques and we obtain as a by
product some regularity results about the solution of the variational
inequality. In particular, we realize there that even if the density
of the underlying diffusion has zeros, there is still a solution on
the open set which corresponds to the region of $[0,T]\times \R^d$
where the density is strictly positive.

%%%%%%%%%%%%%%%%%%%%%%%%%%%%%%%%%%
%%%%%%%%%%%%%%%%%%%%%%%%%%%%%%%%%%%%
\section{Preliminaries and notation}
\label{preliminaries}
\noindent
Let $W$ be the classical Wiener  space $C([0,T],\R^n)$  with  the Wiener
measure $\mu$. The
corresponding Cameron-Martin space is denoted by $H$. Recall that the
injection $H\hookrightarrow W$ is compact and its adjoint is the
natural injection $W^\star\hookrightarrow H^\star\subset
L^2(\mu)$. % A subspace $F$ of $H$ is called regular if the
% corresponding orthogonal projection
% has a continuous extension to $W$, denoted again  by the same letter.
% It is well-known that there exists an increasing sequence of regular
% subspaces $(F_n,n\geq 1)$, called total,  such that $\cup_nF_n$ is
% dense in $H$ and in $W$. Let $\sigma(\pi_{F_n})${\footnote{For the notational
%   simplicity, in the sequel we shall denote  it by  $\pi_{n}$.}}  be the
% $\sigma$-algebra generated by $\pi_{F_n}$, then  for any  $f\in
% L^p(\mu)$, the martingale  sequence
% $(E[f|\sigma(\pi_{F_n})],n\geq 1)$
% converges to $f$ (strongly if  $p<\infty$) in $L^p(\mu)$. Observe that the function
% $f_n=E[f|\sigma(\pi_{F_n})]$ can be identified with a function on the
% finite dimensional abstract Wiener space $(F_n,\mu_n,F_n)$, where
% $\mu_n=\pi_n\mu$.

Since the translations of $\mu$ with the elements of $H$ induce measures
equivalent to $\mu$, the G\^ateaux  derivative in $H$ direction of the
random variables is a closable operator on $L^p(\mu)$-spaces and  this
closure will be denoted by $\nabla$ cf.,  for example
\cite{F-P},\cite{ASU, ASU-1}. The corresponding Sobolev spaces
(the equivalence classes) of the  real random variables
will be denoted as $\DD_{p,k}$, where $k\in \NN$ is the order of
differentiability and $p>1$ is the order of integrability. If the
random variables are with values in some separable Hilbert space, say
$\Phi$, then we shall define similarly the corresponding Sobolev
spaces and they are denoted as $\DD_{p,k}(\Phi)$, $p>1,\,k\in
\NN$. Since $\nabla:\DD_{p,k}\to\DD_{p,k-1}(H)$ is a continuous and
linear operator its adjoint is a well-defined operator which we
represent by $\delta$.  $\delta$ coincides with the It\^o
integral of the Lebesgue density of the adapted elements of
$\DD_{p,k}(H)$ (cf.\cite{ASU,ASU-1}).

For any $t\geq 0$ and measurable $f:W\to \reals_+$, we note by
$$
P_tf(x)=\int_Wf\left(e^{-t}x+\sqrt{1-e^{-2t}}y\right)\mu(dy)\,,
$$
it is well-known that $(P_t,t\in \reals_+)$ is a hypercontractive
semigroup on $L^p(\mu),p>1$,  which is called the Ornstein-Uhlenbeck
semigroup (cf.\cite{F-P,ASU,ASU-1}). Its infinitesimal generator is denoted
by $-\calL$ and we call $\calL$ the Ornstein-Uhlenbeck operator
(sometimes called the number operator by the physicists). The
norms defined by
\begin{equation}
\label{norm}
\|\phi\|_{p,k}=\|(I+\calL)^{k/2}\phi\|_{L^p(\mu)}
\end{equation}
are equivalent to the norms defined by the iterates of the  Sobolev
derivative $\nabla$. This observation permits us to identify the duals
of the space $\DD_{p,k}(\Phi);p>1,\,k\in\NN$ by $\DD_{q,-k}(\Phi')$,
with $q^{-1}=1-p^{-1}$,
where the latter  space is defined by replacing $k$ in (\ref{norm}) by
$-k$, this gives us the distribution spaces on the Wiener space $W$
(in fact we can take as $k$ any real number). An easy calculation
shows that, formally, $\delta\circ \nabla=\calL$, and this permits us
to extend the  divergence and the derivative  operators to the
distributions as linear,  continuous operators. In fact
$\delta:\DD_{q,k}(H\otimes \Phi)\to \DD_{q,k-1}(\Phi)$ and
$\nabla:\DD_{q,k}(\Phi)\to\DD_{q,k-1}(H\otimes \Phi)$ continuously, for
any $q>1$ and $k\in \reals$, where $H\otimes \Phi$ denotes the
completed Hilbert-Schmidt tensor product (cf., for instance
\cite{Mal,ASU,ASU-1}). We shall denote by $\DD(\Phi)$ and $\DD'(\Phi)$
respectively the sets
$$
\DD(\Phi)=\bigcap_{p>1,k\in \N}\DD_{p,k}(\Phi)\,,
$$
and 
$$
\DD'(\Phi)=\bigcup_{p>1,k\in \N}\DD_{p,-k}(\Phi)\,,
$$
where the former is equipped with the projective and the latter is
equipped with the inductive limit topologies. A map $F\in\DD(\R^d)$ is
called nondegenerate if $\det\gamma\in\cap_pL^p(\mu)$, where $\gamma$ is
the inverse of the matrix $((\nabla F_i,\nabla F_j)_H,i,j\leq d)$ and
$(\cdot,\cdot)_H$ denotes the scalar product in $H$. For
such a map, it is well-known that (\cite{Mal,ASU,ASU-1}) the map $f\to f\circ F$ from
$\calS(\R^d)\to \DD$ has a linear, continuous extension to
$\calS'(\R^d)\to \DD'$, where $\calS(\R^d)$ and $\calS'(\R^d)$ denote the space of
rapidly decreasing functions and tempered distributions on $\R^d$,
respectively. In fact, due to the ``polynomially  increasing''
character of the tempered distributions, the range of this extension
is much smaller than $\DD'$, in fact it is included in
$$
\tilde{\DD}'=\bigcap_{p>1}\bigcup_{k\in \NN}\DD_{p,-k}\,.
$$
This notion has been extended in \cite{ASU-0} and used
to give an extension of the It\^o formula as follows: 
\begin{theorem}
\label{ito-thm}
Assume that $(X_t,t\in [0,T])$ is an $\R^d$-valued non-degenerate  It\^o process with
the decomposition 
$$
dX_t=b_tdt+\sigma_t dW_t
$$
where $b\in\DD^a(L^2([0,T])\otimes \R^d)$ and
$\sigma\in\DD^a(L^2([0,T])\otimes \R^d\otimes \R^n)$, where the upper  index
$^a$ means adapted to the Brownian filtration. Assume further that 
\begin{equation}
\label{non-deg-cond}
\int_\eps^1(\det\gamma_s)^pds\in L^1(\mu)\,,
\end{equation}
for any $p>1$, where $\gamma_s$ is the inverse of the matrix $((\nabla
X^i_t,\nabla X^j_t)_H;\,i,j\leq d)$. Then, for any $T\in \calS'(\R^d)$ and $0<s<t\leq 1$, we
have 
$$
T(X_t)-T(X_s)=\int_s^t A_u T(X_u)du+\int_s^t(\partial
T(X_u),\sigma_udW_u)\,,
$$
where $A_u=\frac{1}{2}a_{i,j}(u)\partial_{i,j}+b_i(u)\partial_i$, the
first integral is a Bochner integral in $\tilde{\DD}'$ and the second one is
the extended divergence operator explained above.  
\end{theorem}
\begin{remarkk}
The conditions under which the hypothesis (\ref{non-deg-cond}) holds
are extremely well-studied in the literature, cf. \cite{K-S-1,K-S-2}.
\end{remarkk}
\begin{remarkk}
The divergence operator acts
as an isomorphism between the spaces $\DD_{p,k}^a(H)$ and $\DD_{p,k}$
for any $p>1,\,k\in\R$, cf. \cite{ASU-00}.
\end{remarkk}
\begin{remarkk}
We can extend the above result easily to the case where $t\to T_t$ is
a continuous  map of finite total variation  from $[0,T]$ to $\calS'(\R^d)$ in the sense that,
 the mapping
$t\to \langle T_t,g\rangle$ is of finite total variation on $[0,T]$ for any
$g\in\calS(\R^d)$. In fact, the kernel theorem of A. Grothendieck
implies that $T_t$ can be represented as 
$$
T_t=\sum_{i=1}^\infty \la_i\alpha_i(t)\,F_i\,,
$$
where $(\la_i)\in l^1$, $(\alpha_i)$ is bounded in the total variation
norm and $(F_i)$ is bounded in $\calS'(\R^d)$. Using this
decomposition, it is straightforward to show that 
$$
T(t,X_t)-T(s,X_s)=\int_s^t A_u T(u,X_u)du+ \int_s^tT(du,X_u) +\int_s^t(\partial
T(u,X_u),\sigma_udW_u)\,.
$$
where the second integral is defined as 
$$
\int_s^tT(du,X_u)=\sum_{i=1}^\infty \la_i\int_s^t\,F_i(X_u)
d\alpha_i(u)
$$
and the right hand side is independent of any particular  representation of $T_t$.
integrals are concentrated in $\DD'$
\end{remarkk}

\noindent
We can prove easily  the following result using the technique
described in \cite{ASU-0,ASU-1}: 
\begin{theorem}
\label{ito-1}
Assume that $(l_t,t\in [0,1])$ is an It\^o process 
$$
dl_t=m_tdt+\sum_iz^i_tdW^i_t\,,
$$
with $m,z^i\in \DD^a(L^2[0,T])$, then we have 
\beaa
l_tT(t,X_t)-l_s T(s,X_s)&=&\int_s^t l_uA_u T(u,X_u) du +\int_s^t l_u T(du,X_u)\\
&&+\int_s^tl_u(\partial
T(u,X_u),\sigma_udW_u)+\int_s^tT(u,X_u)m_u du\\
&&+\int_s^tT(u,X_u)\sum_iz^i_udW_u^i+\int_s^t(\partial
T(u,X_u),\sigma_u z_u)du
\eeaa
where $z=(z^1,\ldots,z^n)$.
\end{theorem}

An important feature of the distributions on the Wiener space is the
notion of positivity: we say that $S\in\DD'$ is positive if for any positive
$\varphi\in\DD$, we have $S(\varphi)=\langle S,\varphi\rangle\geq
0$. An important result about the positive distributions is the
following (cf. \cite{A-M,Mal,Sug,ASU,ASU-1}):
\begin{theorem}
Assume that $S$ is a positive distribution in $\DD'$, then there
exists a positive Radon measure $\nu_S$ on $W$ such that 
$$
\langle S,\varphi\rangle=\int_W \varphi d\nu_S\,,
$$
for any $\varphi\in \DD\cap C_b(W)$. In particular, if  a sequence $(S_n)$
of positive distributions converge to $S$ weakly in $\DD'$, then
$(\nu_{S_n})$ converges to $\nu_S$ in the weak topology of measures.
\end{theorem}
\begin{remarkk}
In fact we can write, for any $\varphi\in \DD$
$$
\langle S,\varphi\rangle=\int_W\tilde{\varphi} d\nu_S\,,
$$
where $\tilde{\varphi}$ denotes a redefinition of $\varphi$ which is
constructed  using the
capacities associated to the scale of Sobolev spaces
$(\DD_{p,k},\,p>1,k\in\N)$, cf. \cite{F-P}.
\end{remarkk}

\section{Uniqueness of  the  solution of parabolic variational inequality}
Assume that $(X^s_t(x), \,0\leq s\leq t\leq T)$ is a diffusion process
governed by an $\R^n$-valued Wiener process $(W_t,t\in [0,T])$.
We assume that the diffusion has  smooth, bounded  drift and diffusion coefficients $b(t,x),\sigma(t,x)$
defined on $[0,T]\times \R^d$, with values in $\R^d$ and $\R^d\otimes
\R^n$ respectively and we denote by $A_t$ its infinitesimal
generator. We shall assume that $X_t^s$ is nondegenerate for any
$0\leq s<t\leq T$, $\partial/\partial t+A_t$ is  hypoelliptic
and 
$$
\int_{s+\eps}^t (\det\gamma^s_v)^pdv\in L^1(\mu)
$$
for any $0<s<t \leq T$ and $\eps>0$, where $\gamma_v$ is the inverse of
the matrix $((\nabla X^{s,i}_v,\nabla X^{s,j}_v)_H:\,i,j\leq d)$.

Suppose that $f\in C_b(\R^d)$ and we shall study the
following partial differential inequality whose solution will be
denoted by $u(t,x)$:
\begin{theorem}
\label{ineq-th-1}
Assume that $u\in C_b([0,T]\times \R^d)$ such that $t\to
\langle u(t,\cdot),g\rangle $ is of finite total variation on $[0,T]$
for any $g\in\calS(\R^d)$ and that it  satisfies the following properties:
\begin{eqnarray}
\frac{\partial u}{\partial t}+A_tu-ru&\leq& 0,\,\,u\geq
f\,{\mbox{in}}\,[0,T]\times \R^d\label{rel-1}\\
(\frac{\partial u}{\partial t}+A_tu-ru)(f-u)&=& 0,\,\label{rel-2}\\
U(T,x)&=&f(x)\,,\label{rel-3}
\end{eqnarray}
where all the derivatives are taken in the sense of distributions, in
particular the derivative w.r. to $t$ is taken using the
$C^\infty$-functions of compact support in $(0,T)$.
Then
$$
u(t,x)=\sup_{\tau\in\calZ_{t,T}}E\left[f(X^{t}_\tau(x))\exp-\int_t^\tau r(s,X_s^t(x))ds\right]\,,
$$
where $\calZ_{t,T}$ denotes the set of all the stopping times with
values in $[t,T]$ and $r$ is a smooth function on $[0,T]\times \R^d$.
\end{theorem}
\nproof
We shall prove the case $t=0$. Let us denote by $l$ the process
defined as $l_t=\exp-\int_0^tr(s,X_s)ds$. From Theorem \ref{ito-1}, we
have, for any $\eps>0$, 
\begin{equation}
\label{ito-dev}
l_tu(t,X_t)-l_\eps u(\eps,X_\eps)-\int_\eps^t l_s (A_su(s,X_s)ds-(ru)(s,X_s)ds+u(ds,X_s))=M^\eps_t
\end{equation}
where  $M^\eps_t$ is a
$\DD'$-valued martingale difference, i.e., denoting by
$E[\cdot|\calF_s]$ the extension of the conditional expectation
operator to $\DD'$ {\footnote{Such an extension is licit since the
    conditional expectation
    operator commutes with the Ornstein-Uhlenbeck semigroup.}}, we
have $E[M^\eps_t|\calF_s]=M^\eps_s$ for any $\eps\leq s\leq t$. Note
also that $K_tu=\frac{\partial u}{\partial t}+A_tu-ru)\leq 0$ hence its composition with $X_t$ is a negative
measure and this implies that the integral at the l.h.s. of
(\ref{ito-dev}) is a negative distribution on the Wiener space. Consequently we have 
\begin{equation}
\label{eqn-1}
M^\eps_t\geq l_tu(t,X_t)-l_\eps u(\eps,X_\eps)
\end{equation}
in $\DD'$. For $\alpha>0$, let $P_\alpha$ be the Ornstein-Uhlenbeck
semigroup and define $M^{\alpha,\eps}_t$ as
$$
M^{\alpha,\eps}_t=P_\alpha M^\eps_t\,.
$$
Then $(M^{\alpha,\eps}_t,t\geq \eps)$ is a continuous martingale (in
the ordinary sense). From the inequality (\ref{eqn-1}), we have, for
any $\tau\in\calZ_{\eps,T}$, 
$$
M^{\alpha,\eps}_\tau\geq P_\alpha(l_tu(t,X_t)-l_\eps
u(\eps,X_\eps))|_{t=\tau}\,.
$$
Taking the expectation of both sides, we get
$$
E[l_\eps u(\eps,X_\eps)]\geq E[ P_\alpha(l_tu(t,X_t)|_{t=\tau}]
$$
for any $\alpha>0$, hence we also have 
$$
E[l_\eps u(\eps,X_\eps)]\geq E[ l_\tau u(\tau,X_\tau)]
$$
for any $\eps>0$ which is arbitrary and finally  we obtain
$$
u(0,x)\geq E[ l_\tau u(\tau,X_\tau)]
$$
for any $\tau\in \calZ_{0,T}$.

\noindent
To show the reverse inequality let $D=\{(s,x):\,u(s,x)\neq f(x)\}$ and
define 
$$
\tau_x=\inf(s:\,(s,X^{0,x}_s)\in D^c)\,.
$$
 Since $K_t$ is
hypoelliptic, and since  $K_tu=0$ on the set $D$, $u$ is smooth in
$D$. If $\mu\{\tau_x=0\}=1$, from the continuity of $u$, we have 
$$
u(0,x)=f(x)=E[l_{\tau_x}u(\tau_x,X^{0,x}_{\tau_x})]\,,
$$
hence the supremum is attained in this case. If $\mu\{\tau_x\neq
0\}>0$, then from the $0-1$-law $\mu\{\tau_x\neq 0\}=1$ and $\tau_x$
is predictable. Let $(\tau_n,n\geq 1)$ a sequence of stopping times
announcing $\tau_x$. From the classical It\^o formula, we have 
$$
l_{\tau_n}u(\tau_n,X_{\tau_n})-u(0,x)=\int_0^{\tau_n}l_s(\sigma^\star\partial
u)(s,X_s)\cdot dW_s\,.
$$
By the hypothesis the l.h.s. is uniformly integrable with respect to
$n\in\NN$, consequently we obtain 
$$
u(0,x)=\lim_nE[l_{\tau_n}u(\tau_n,X_{\tau_n})]=E[l_{\tau}u(\tau,X_{\tau})]
$$
hence $\tau_x$ realizes the supremum.
\nqed

\noindent
In the homogeneous case the finite variation property of the solution
follows directly from the quasi-variational inequality:

\begin{theorem}
Suppose that the infinitesimal generator  $A_t$ of the process $(X_t)$ is
independent of $t\in [0,T]$ and denote it by $A$. In other words  the process is homogeneous in
time. Assume that $u\in C_b([0,T]\times \R^d)$   satisfies the following properties:
\begin{eqnarray}
\frac{\partial u}{\partial t}+Au-ru&\leq& 0,\,\,u\geq
f\,{\mbox{in}}\,[0,T]\times \R^d\label{rel-11}\\
(\frac{\partial u}{\partial t}+Au-ru)(f-u)&=& 0\,\,{\mbox{in}}\,[0,T]\times \R^d\label{rel-12}\\
U(T,x)&=&f(x)\,.\label{rel-13}
\end{eqnarray}
Then
$$
u(t,x)=\sup_{\tau\in\calZ_{t,T}}E\left[f(X^{t}_\tau(x))\exp-\int_t^\tau r(s,X_s^t(x))ds\right]\,,
$$
where $\calZ_{t,T}$ denotes the set of all the stopping times with
values in $[t,T]$ and $r$ is a smooth function on $[0,T]\times \R^d$.
\end{theorem}
\begin{remarkk}
The relations (\ref{rel-11}) and (\ref{rel-12}) are to be understood
in the weak sense. This means that for any $g$ a $C^\infty$ function
of support in $(0,T)$ and $\gamma\in\calS(\R^d)$, both of which are
positive, we have 
$$
\left\langle \frac{\partial u}{\partial t}+Au-ru,g\otimes\gamma\right\rangle \leq 0
$$
and 
$$
\left\langle (\frac{\partial u}{\partial t}+Au-ru)(f-u),g\otimes\gamma\right\rangle= 0\,.
$$
\end{remarkk}
\nproof
As in the proof of the preceding theorem, we shall prove the equality
for $t=0$, then  the general case follows easily.
Let $\rho_\delta$ be a mollifier on $\R$ and let $\eta_\eps$ be a
family of positive  smooth functions on $(0,T)$, equal to unity on the
interval $[\eps,T-\eps]$, converging to the indicator function of
$[0,T]$ pointwise. Define $u^{\delta,\eps}$ as 
$$
u^{\delta,\eps}=\rho_\delta\star(\eta_\eps u)\,.
$$
From the hypothesis the distribution $\nu$ defined by
$$
\nu=\frac{\partial u}{\partial t}+Au-ru
$$
is a negative measure on $(0,T)\times \R^d$. A simple calculation
gives
$$
\frac{\partial u^{\delta,\eps}}{\partial
  t}+Au^{\delta,\eps}-ru^{\delta,\eps}=\rho_\delta\star(\eta'_\eps
u)+\rho_\delta\star(\eta_\eps\nu)+\rho_\delta\star(\eta_\eps r
u)-ru^{\delta,\eps}\,.
$$
As in the preceding theorem, we have from Theorem \ref{ito-1}
$$
l_tu^{\delta,\eps}(t,X_t)-l_au^{\delta,\eps}(a,X_a)-\int_a^tl_sK_su^{\delta,\eps}(s,X_s)ds=M_t^{\delta,\eps,a}\,,
$$
where $M^{\delta,\eps,a}$ is a $\DD'$-martingale difference. Since
$\nu$ is a negative measure, we get the following inequality in
$\DD'$:
$$
M_t^{\delta,\eps,a}\geq l_tu^{\delta,\eps}(t,X_t)-l_au^{\delta,\eps}(a,X_a)-\int_a^t[\rho_\delta\star(\eta'_\eps
u)+\rho_\delta\star(\eta_\eps r u)-ru^{\delta,\eps}](s,X_s)ds
$$
Let now $(P_\alpha,\alpha\geq 0)$ be the Ornstein-Uhlenbeck
semigroup. Then $(P_\alpha M_t^{\delta,\eps,a},\,a\leq t\leq T)$ is a
real valued martingale difference, consequently, we have 
\beaa
0&=&E[(P_\alpha M_t^{\delta,\eps,a})_{t=\tau}]\\
&\geq&
E\left[P_\alpha\left(l_t
    u^{\delta,\eps}(t,X_t)-l_au^{\delta,\eps}(a,X_a)-
\int_a^t[\rho_\delta\star(\eta'_\eps
u)+\rho_\delta\star(\eta_\eps r
u)-ru^{\delta,\eps}](s,X_s)ds\right)_{t=\tau}\right]\,,
\eeaa
for any stopping time $\tau$ with values in $[\eps,T-\eps]$. By
letting $\alpha\to 0$, we get by continuity
$$
0\geq
E\left[l_\tau
    u^{\delta,\eps}(\tau,X_\tau)-l_a u^{\delta,\eps}(a,X_a)-
\int_a^\tau[\rho_\delta\star(\eta'_\eps
u)+\rho_\delta\star(\eta_\eps r
u)-ru^{\delta,\eps}](s,X_s)ds\right]\,.
$$
 Let us
choose $a>0$ and let then $\eps,\delta\to 0$. Note that
$\eta'_\eps \to \delta_0-\delta_T$ (i.e., the Dirac measures at $0$
and at $T$), by the choice of $a$ and by the weak convergence of
measures and by the dominated convergence theorem, we obtain
$$
\lim_{\eps,\delta\to 0}E\int_a^\tau \left(\rho_\delta\star(\eta'_\eps
u)\right)(s,X_s)ds=0\,.
$$
Again from the dominated convergence theorem we have 
$$
\lim_{\eps,\delta\to 0}E\int_a^\tau [\rho_\delta\star(\eta_\eps r
u)-ru^{\delta,\eps}](s,X_s)ds=0\,.
$$
Consequently
$$
E[l_a u(a,X_a)]\geq E[l_\tau u(\tau,X_\tau)]\geq E[l_\tau f(X_\tau)]\,,
$$
for any stopping time $\tau$ with values in $[a,T-a]$, since $a>0$ is
arbitrary, the same inequality holds also for any stopping time with
values in $[0,T]$; hence 
$$
 u(0,x)]\geq E[l_\tau f(X_\tau)]
$$
for any stopping time $\tau\in \calZ_{0,T}$ and we obtain the first
inequality:
$$
u(0,x)]\geq\sup_{\tau\in\calZ_{0,T}} E[l_\tau f(X_\tau)]\,.
$$
The proof of the reverse inequality is exactly the same that of
Theorem \ref{ineq-th-1} due to the
hypoellipticity  hypothesis.
\nqed

\section{Existence of the solutions}

\noindent
In this section, under the hypothesis of the preceding section,  we shall prove that the function defined  by the
Snell envelope (cf. \cite{Karoui}) of the American option satisfies the variational
inequality (\ref{rel-11}) and the equality (\ref{rel-12}).
We start with a lemma:
\begin{lemma}
\label{lemma-1}
Assume that $Z=(Z_t,t\in[0,T])$ is a uniformly integrable, real-valued  martingale on
the Wiener space. Let $Z^\kappa=(Z^\kappa_t,t\in[0,T])$ be defined as
$Z_t^\kappa=P_\kappa Z_t$, where $P_\kappa$ is the Ornstein-Uhlenbeck
semigroup at the instant $\kappa>0$. Then $(Z_t^\kappa,t\in[0,T])$  is
a uniformly integrable martingale with
\begin{equation}
\label{Davis}
E[\langle Z^\kappa, Z^\kappa\rangle_T^{1/2}]\leq c\,E[\langle Z,
Z\rangle_T^{1/2}]\,,
\end{equation}
where $c$ is a constant independent of $Z$ and $\kappa$. In
particular, if $Z$ has the representation 
$$
Z_T=\int_0^T(m_s,dW_s)\,,
$$
with $m\in L^1(\mu,H)$ optional, then 
$$
P_\kappa Z_T=\int_0^Te^{-\kappa}(P_\kappa m_s,dW_s)\,.
$$
\end{lemma}
\nproof
From Davis'  inequality (cf. \cite{Mey}), we have 
\beaa
E[\langle Z^\kappa, Z^\kappa\rangle_T^{1/2}]&\leq&c_1 E[\sup_{t\in [0,T]}|Z_t^\kappa|]\\
&\leq&c_1 E[P_\kappa(\sup_{t\in [0,T]}|Z_t|)]\\
&=&c_1 E[\sup_{t\in [0,T]}|Z_t|]\\
&\leq& c\,E[\langle Z,Z\rangle_T^{1/2}]\,.
\eeaa
The second part is obvious from the inequality (\ref{Davis}).
\nqed

\begin{theorem}
\label{existence-thm}
Assume that $(X_t^s)$ is a hypoelliptic diffusion such that, for any
$\eps>0$, 
$$
\int_{s+\eps}^T (\det\gamma_v^s)^p dv\in L^1(\mu)
$$
for any $p>1$. Let $p(s,t;x,y),\,s<t,\,x,y\in \R^d$ be the density of
the law of $X^s_t(x)$ and denote by $S_{0,z}$ the open set
$$
S_{0,z}=\{(s,y)\in (0,T)\times \R^d:\,s>0,\,p(0,s;z,y)>0\}\,.
$$
Then, for any $z\in \R^d$,  $u$ is a solution of the variational inequality
(\ref{rel-11},\ref{rel-12},\ref{rel-13}) in $\calD'(S_{0,z})$. If
$S_{0,z}=(0,T)\times \R^d$ for any $z\in \R^d$, then $u$ is a solution of the variational inequality
(\ref{rel-11},\ref{rel-12},\ref{rel-13}) in $\calD'(0,T)\otimes \calD'(\R^d)$.
\end{theorem}
\nproof
From the optimal stopping results, we know that $u$ is a bounded,
continuous function and  $t\to u(t,x)$ is monotone, decreasing
(cf.\cite{Karoui}). Moreover
$$
u(t,X_t)l_t-u(0,x)=M_t+B_t
$$
is a supermartingale where $X_t=X^0_t(x)$ and  we denoted by $M$ its
martingale part and by $B$ its continuous, decreasing process part. In
particular $dB\times d\mu$ defines a negative  measure $\gamma$  on
$[0,T]\times C([0,T],\R^d)$. We can write $u(ds,x)$ as the sum
$u_{ac}(s,x)ds+u_{sing}(ds,x)$ where $u_{ac}$ is defined as the absolutely
continuous part of $u$ and $u_{sing}$ is the singular part.
 We have, from the extended It\^o formula, 
\beaa
u(t,X_t)l_t-u(\eps,X_\eps)&=&\int_\eps^t
\left(A_su-ru+u_{ac})(s,X_s)\right)ds+u_{sing}(ds,X_s)+\int_\eps^t((\sigma\partial
u)(s,X_s),dW_s)\\
&=&M^\eps_t+B^\eps_t\,,
\eeaa
hence regularizing both parts by the Ornstein-Uhlenbeck semigroup,
from Lemma \ref{lemma-1},  we get  
\beaa
B^\eps_t&=&\int_\eps^t \left(A_su-ru+u_{ac})(s,X_s)\right)ds+u_{sing}(ds,X_s)\\
M^\eps_t&=&\int_\eps^t((\sigma\partial u)(s,X_s),dW_s)\,.
\eeaa
Consequently, for any  $\alpha\in \calD(0,T)$ and $\phi\in \DD$
\beaa
E\left[\phi\int_0^T\alpha(s) dB_s\right]&=&\int\alpha\otimes \phi d\gamma\\
&=&\int_{(0,T)}\alpha(s)\langle
\left(A_su-ru+u_{ac})(s,X_s)\right)ds+u_{sing}(ds,X_s),\phi\rangle 
\eeaa
and this quantity is negative for any $\alpha\in \calD_+(0,T)$ and
$\phi\in \DD_+$. Let now $0\leq g\in\calS(\R^d)$ and assume that
$(t_i,\,i\leq m)$ is a partition of $[0,T]$. Define $\xi_m$ as 
$$
\xi_m(t,w)=\sum_i 1_{[t_i,t_{i+1}]}(t)\,g(X_{t_i})\,.
$$
Then it is immediate from the hypothesis about the diffusion process
$(X_t)$ that $(\xi_m,\,m\geq 1)$ converges to
$(g(X_s)1_{[0,T]}(s),\,s\in [0,T])$ in $\DD(L^p([0,T]))$ for any
$p\geq 1$ and $(\xi_m(s,\cdot)$ converges to $g(X_s)$ in $\DD$  for
any fixed $s\in [0,T]$ as the partition pace tends to zero. Let us
represent $u(ds,\cdot)$, using the kernel theorem (cf.\cite{Gro,Sch}), as 
$$
u(ds,\cdot)=\sum_{k=1}^\infty \la_i T_k\otimes \alpha_k\,,
$$
where $(\la_k)\in l^1$, $(T_k)\subset \calS'(\R^d)$ is bounded and
$(\alpha_k)$ is a sequence of measures on $[0,T]$, bounded in total
variation norm. It follows then
$$
u(ds,X_s)=\sum_{k=1}^\infty \la_i T_k(X_s) \alpha_k(ds)
$$
and this some is convergent in $V([0,T])\tilde{\otimes}\DD_{p,-k}$ for
some $k\in\N$ and $p>1$, in the projective topology, where $V([0,T])$ denotes the
Banach space of measures on $[0,T]$ under the total variation norm.
Since 
$$
\sup_{s\in [0,T]}\|\xi_m(s,X_s)\|_{p,l}\leq \sup_{s\in [0,T]}\|g\circ
X_s\|_{p,l} \,,
$$
uniformly in $m\in \NN$, for any $p,\,l$ and since 
$$
\|\xi_m(s,\cdot)-\xi_n(s,\cdot)\|_{p,l}\to 0
$$
as $m,n\to\infty$ for any $p,\,l$ and $s\in [0,T]$,  we obtain 
$$
\lim_{m\to\infty}\int_{(0,T)}\delta(s)\langle
\xi_m(s,\cdot)-g(X_s),u(ds,X_s)\rangle=0
$$
for any $\delta\in\calD(0,T)$ from the dominated convergence theorem.
The above relation implies in particular
 that we  have 
$$
\int_{(0,T)}\alpha(s)\langle (A_su-ru+u_{ac})(s,\cdot),p_{0,s}g\rangle
ds+\int_{(0,T)}\alpha(s)\langle u_{sing}(ds,\cdot),g\,p_{0,s}\rangle\leq 0\,,
$$
 with smooth, positive $\alpha$ and $g$, where the brackets in the integral correspond to the duality between
$\calD(\R^d)$ and $\calD'(\R^d)$. For the functions of support in
$(0,T)$, we can replace the term 
$$
\int_{(0,T)} \alpha(s)\left\langle u(ds,X_s),g(X_s)\right\rangle
$$
by 
$$
\int_{(0,T)} \alpha(s)\langle\frac{\partial}{\partial s} u(s,X_s),g(X_s)\rangle
$$
where $\partial/\partial s$ denotes the derivative in $\calD'(0,T)$.
 Since $\alpha$ and $g$ are
arbitrary, we obtain the inequality (\ref{rel-11}) in
$\calD'(S_{0,x})$. If $S_{0,x}=(0,T)\times \R^d$, then we have the
inequality in the sense of distributions on $(0,T)\times \R^d$.

To complete the proof, let $D$ be the set defined as
$$
D=\{(s,x)\in (0,T)\times \R^d: \,u(s,x)=f(x)\}\,.
$$
Then we have 
$$
\int_0^T 1_{D^c}(s,X_s) dB_s=0
$$
almost surely (cf.\cite{Karoui}). Let $C=-B$, then for any smooth
function $\eta\in \calD(0,T)\otimes \calS(\R^d)$ such that $\eta\leq
1_{D^c}$, we have 
\beaa
0&=&E\int_0^T1_{D^c}(s,X_s)dC_s\\
&\geq&E\int_0^T\eta(s,X_s)dC_s\\
&=&-\int_{(0,T)}\langle (A_s u-r u)(s,X_s)ds+u(ds,X_s),\eta(s,X_s)\rangle\\
&\geq&0\,,
\eeaa
where, the second equality follows from the estimates above. Hence 
$$
A_su-ru+\frac{\partial}{\partial s}u=0
$$
as a distribution on the set $S_{0,x}\cap D^c$, by the
hypoellipticity, the equality is everywhere on this set. If
$S_{0,x}=(0,T)\times \R^d$, then we obtain the relation
(\ref{rel-12}).
\nqed
\begin{remarkk}
From the general theory, we can express the martingale part of
$(l_tu(t,X_t),\,t\in [0,T])$ as 
$$
\int_0^T(H_s,dW_s)
$$
where $H$ is an adapted process which is locally integrable. On the
other hand we have 
$$
M^\eps_t=\int_\eps^t (\sigma(s,X_s)\partial u(s,X_s),dW_s)
$$
where the r.h.s. is to be interpreted in a negatively indexed Sobolev
space on the Wiener space. Using Lemma \ref{lemma-1}, we obtain the
identity 
$$
H_s=\sigma(s,X_s)\partial u(s,X_s)
$$
$ds\times d\mu$-a.s., in particular we have 
$$
E\left[\left(\int_0^T|\sigma(s,X_s)\partial u(s,X_s)|^2ds\right)^{1/2}\right]<\infty\,.
$$
\end{remarkk}

%\newpage

\vspace{2cm}

{\footnotesize{\bf{
\noindent
A.S. \"Ust\"unel, Telecom-Paristech (formerly ENST),
 Dept. Infres,\\
46, rue Barrault, 75013 Paris, France\\
email: ustunel@telecom-paristech.fr}
}}

\end{document}